\documentclass[10pt]{amsart}
\usepackage[cp1251]{inputenc}
\usepackage[english]{babel}
\usepackage{amsmath}
\usepackage{amssymb}
\usepackage{amsfonts}

\usepackage[linktocpage=true, colorlinks=true, linkcolor=blue, citecolor=blue, urlcolor=blue]{hyperref}

\begin{document}
\renewcommand{\refname}{References}

\thispagestyle{empty}

\title[Expansion of Iterated Stratonovich Stochastic Integrals of 
Multiplicity 2]
{Expansion of Iterated Stratonovich Stochastic Integrals of Multiplicity 2.
Combined Approach Based on Generalized Multiple and Iterated Fourier
Series}
\author[D.F. Kuznetsov]{Dmitriy F. Kuznetsov}
\address{Dmitriy Feliksovich Kuznetsov
\newline\hphantom{iii} Peter the Great Saint-Petersburg Polytechnic University,
\newline\hphantom{iii} Polytechnicheskaya ul., 29,
\newline\hphantom{iii} 195251, Saint-Petersburg, Russia}%
\email{sde\_kuznetsov@inbox.ru}

\vspace{5mm}

\maketitle {\small
\begin{quote}
\noindent{\sc Abstract.} 
The article is devoted to the expansion of iterated
Stratonovich stochastic integrals of multiplicity 2
on the base of the combined approach of generalized multiple and iterated
Fourier series.
We consider two different parts of the expansion of iterated
Stra\-to\-no\-vich stochastic integrals.
The mean-square convergence of the first part is proved on the base
of generalized multiple Fourier series converging 
in the sense of norm in Hilbert space 
$L_2([t, T]^2).$ The mean-square convergence
of the second part is proved on the base of 
generalized iterated (double) Fourier
series converging pointwise. At that, we prove the iterated 
limit transition
for the second part of the expansion on the base of 
Lebesgue's Dominated Convergence Theorem.
The results of the article can be applied to the numerical integration 
of Ito stochastic differential equations.

\vspace{5mm}


\setlength{\baselineskip}{2.0em}

\tableofcontents

\setlength{\baselineskip}{1.2em}


\medskip
\end{quote}
}

\section{Introduction}

\vspace{5mm}

Let $(\Omega,$ ${\rm F},$ ${\sf P})$ be a complete probability space, let 
$\{{\rm F}_t, t\in[0,T]\}$ be a nondecreasing right-continous family of 
$\sigma$-algebras of ${\rm F},$
and let ${\bf f}_t$ be a standard $m$-dimensional 
Wiener stochastic process, which is
${\rm F}_t$-measurable for any $t\in[0, T].$ We assume that the components
${\bf f}_{t}^{(i)}$ $(i=1,\ldots,m)$ of this process are independent.

Let us consider the following collections of iterated
Stratonovich and Ito stochastic integrals

\vspace{-1mm}
\begin{equation}
\label{str}
J^{*}[\psi^{(2)}]_{T,t}=
{\int\limits_t^{*}}^T\psi_2(t_2){\int\limits_t^{*}}^{t_2}
\psi_1(t_1) d{\bf w}_{t_1}^{(i_1)}d{\bf w}_{t_2}^{(i_2)},
\end{equation}

\begin{equation}
\label{ito}
J[\psi^{(2)}]_{T,t}=\int\limits_t^T\psi_2(t_2)\int\limits_t^{t_{2}}
\psi_1(t_1) d{\bf w}_{t_1}^{(i_1)}
d{\bf w}_{t_2}^{(i_2)},
\end{equation}

\vspace{2mm}
\noindent
where every $\psi_l(\tau)$ $(l=1,\ 2)$ is a nonrandom function 
at the interval $[t,T],$ ${\bf w}_{\tau}^{(i)}={\bf f}_{\tau}^{(i)}$
for $i=1,\ldots,m$ and
${\bf w}_{\tau}^{(0)}=\tau,$\ \
$i_1,\ldots,i_k = 0,\ 1,\ldots,m,$

\vspace{-1mm}
$$
\int\limits^{*}\ \hbox{and}\ \int\limits
$$ 

\vspace{2mm}
\noindent
denote 
Stratonovich and Ito stochastic integrals,
respectively (in this paper, 
we use the definition of the Stratonovich stochastic integral from \cite{KlPl2}).

Further, we will 
denote as $\{\phi_j(x)\}_{j=0}^{\infty}$ the
complete orthonormal systems of Legendre polynomials and
trigonometric functions 
in the space $L_2([t, T])$.
Also we will pay a special attention on the 
following well-known facts connecting
to these two systems of functions \cite{Gob}.

{\it Suppose that the function $f(x)$ is 
bounded at the interval $[t, T].$ Moreover, its derivative
$f'(x)$ is a continuous function at the interval $[t, T]$ except may be
the finite number of points 
of the finite discontinuity.
Then the Fourier series 

\vspace{-1mm}
$$
\sum\limits_{j=0}^{\infty}
C_j\phi_j(x),\ \ \ C_j=\int\limits_t^{T}f(x)\phi_j(x)dx
$$

\vspace{3mm}
\noindent
converges at any internal point $x$ of 
the interval $[t, T]$ to the value 
$\left(f(x+0)+f(x-0)\right)/2$ and converges  
uniformly to $f(x)$ on any closed interval of continuity 
of the function
$f(x)$ laying inside 
$[t, T]$. At the same time, the Fourier--Legendre series 
converges 
if $x=t$ and $x=T$ to $f(t+0)$ and $f(T-0)$ 
correspondently, and the trigonometric Fourier series converges if   
$x=t$ and $x=T$ to $\left(f(t+0)+f(T-0)\right)/2$ 
in the case of periodic continuation 
of the function $f(x)$.}

\vspace{5mm}

\section{Expansion of Iterated Stratonovich Stochastic Integrals of 
Multiplicity 2}

\vspace{5mm}

The use of generalized  multiple and iterated Fourier series 
by various complete orthonormal systems of 
functions in the space $L_2([t, T])$ 
for the expansion of iterated
Ito and Stratonovich stochastic integrals
is reflected in a number of works of the author
\cite{2013}-\cite{01}. In these papers, several new approaches
to the mean-square approximation of iterated
stochastic integrals were proposed and developed.
One of the mentioned approaches (the so-called combined approach)
for the expansion of iterated
Stratonovich stochastic integrals of multiplicities 1 to 4
based on generalized multiple and iterated Fourier series
has been 
considered in \cite{2013a}. In this article, we consider
the case of second multiplicity of iterated
Stratonovich stochastic integrals. At that, we prove the
mean-square convergence of the expansion 
of iterated
Stratonovich stochastic integrals
using the another
method in comparison with the method from
\cite{2013a}.

\vspace{2mm}

{\bf Theorem~1}\ \cite{2013} (2013) (also see \cite{2017-1xx} (Sect.~2.1.1) and references therein).\
{\it Suppose that 
$\{\phi_j(x)\}_{j=0}^{\infty}$ is a complete orthonormal system of 
Legendre polynomials or trigonometric functions in the space $L_2([t, T]).$
At the same time $\psi_2(\tau)$ is a continuously dif\-ferentiable 
nonrandom function on $[t, T]$ and $\psi_1(\tau)$ is twice 
continuously differentiable nonrandom function on $[t, T]$. 
Then the iterated Stratonovich stochastic integral of the second multiplicity

\vspace{-1mm}
$$
J^{*}[\psi^{(2)}]_{T,t}={\int\limits_t^{*}}^T\psi_2(t_2)
{\int\limits_t^{*}}^{t_2}\psi_1(t_1)d{\bf w}_{t_1}^{(i_1)}
d{\bf w}_{t_2}^{(i_2)}\ \ \ (i_1, i_2=0, 1,\ldots,m)
$$

\vspace{2mm}
\noindent
is expanded into the 
multiple series

\vspace{-1mm}
$$
J^{*}[\psi^{(2)}]_{T,t}=
\hbox{\vtop{\offinterlineskip\halign{
\hfil#\hfil\cr
{\rm l.i.m.}\cr
$\stackrel{}{{}_{p_1, p_2\to \infty}}$\cr
}} }
\sum\limits_{j_1=0}^{p_1}
\sum\limits_{j_2=0}^{p_2}
C_{j_2j_1}\zeta_{j_1}^{(i_1)}\zeta_{j_2}^{(i_2)}
$$

\vspace{3mm}
\noindent
that converges 
in the mean-square sense,
where {\rm l.i.m.} is a limit in the mean-square sense,

\vspace{-1mm}
$$
\zeta_{j}^{(i)}=
\int\limits_t^T \phi_{j}(\tau) d{\bf w}_{\tau}^{(i)}
$$

\vspace{2mm}
\noindent
are independent standard Gaussian random variables
for various
$i$ or $j$ {\rm(}if $i\ne 0${\rm)},

\vspace{-1mm}
$$
C_{j_2 j_1}=
\int\limits_t^T
\psi_2(t_2)\phi_{j_2}(t_2)\int\limits_t^{t_2}
\psi_1(t_1)\phi_{j_1}(t_1)dt_1 dt_2
$$

\vspace{2mm}
\noindent
is the Fourier coefficient.}

\vspace{2mm}

{\bf Remark 1.}\ {\it It should be noted that
Theorem {\rm 1} is proved in 
{\rm \cite{2013} (2013) (}also see {\rm \cite{2017-1xx} (}Sect.~{\rm 2.1.1)} and references therein{\rm )}.
The proof from {\rm \cite{2013}, \cite{2017-1xx} (}Sect.~{\rm 2.1.1)} 
is based on double integration by parts.
Below we consider another proof of Theorem {\rm 1}.}

\vspace{2mm}

{\bf Proof.} Let us consider some
auxiliary lemmas from 
\cite{2013} (also see \cite{2017-1xx} and references therein).
At that, we will consider
the particular case of these lemmas for $k=2.$

Consider the partition $\{\tau_j\}_{j=0}^N$ of the interval $[t,T]$ such that

\vspace{-1mm}
\begin{equation}
\label{1111}
t=\tau_0<\ldots <\tau_N=T,\ \ \ 
\Delta_N=
\hbox{\vtop{\offinterlineskip\halign{
\hfil#\hfil\cr
{\rm max}\cr
$\stackrel{}{{}_{0\le j\le N-1}}$\cr
}} }\Delta\tau_j\to 0\ \ \hbox{if}\ \ N\to \infty,\ \ \ 
\Delta\tau_j=\tau_{j+1}-\tau_j.
\end{equation}

\vspace{3mm}

{\bf Lemma 1}\ \cite{2013} (also see \cite{2017-1xx} and references therein).
{\it Suppose that
every $\psi_l(\tau)$ $(l=1,\ 2)$ is a continuous 
nonrandom function at the interval
$[t, T]$. Then

\vspace{-1mm}
\begin{equation}
\label{30.30}
J[\psi^{(2)}]_{T,t}=
\hbox{\vtop{\offinterlineskip\halign{
\hfil#\hfil\cr
{\rm l.i.m.}\cr
$\stackrel{}{{}_{N\to \infty}}$\cr
}} }
\sum_{j_2=0}^{N-1}\sum_{j_1=0}^{j_{2}-1}
\psi_1(\tau_{j_1})\psi_2(\tau_{j_2})   
\Delta{\bf w}_{\tau_{j_1}}^{(i_1)}\Delta{\bf w}_{\tau_{j_2}}^{(i_2)}\ \ \
\hbox{\rm w. p. 1},
\end{equation}

\vspace{2mm}
\noindent
where $J[\psi^{(2)}]_{T,t}$ is the iterated Ito stochastic integral
{\rm (\ref{ito}),} 
$\Delta{\bf w}_{\tau_{j}}^{(i)}=
{\bf w}_{\tau_{j+1}}^{(i)}-{\bf w}_{\tau_{j}}^{(i)}$
$(i=0, 1,\ldots,m)$,
$\left\{\tau_{j}\right\}_{j=0}^{N}$ is the partition 
of the interval $[t,T]$ satisfying the condition {\rm (\ref{1111});}
hereinafter w.~p.~{\rm 1}  means  with probability {\rm 1}.
}

\vspace{2mm}

Let us define the following
multiple stochastic integral

\vspace{-1mm}
\begin{equation}
\label{30.34}
\hbox{\vtop{\offinterlineskip\halign{
\hfil#\hfil\cr
{\rm l.i.m.}\cr
$\stackrel{}{{}_{N\to \infty}}$\cr
}} }\sum_{j_1,j_2=0}^{N-1}
\Phi\left(\tau_{j_1},\tau_{j_2}\right)
\Delta{\bf w}_{\tau_{j_1}}^{(i_1)}\Delta{\bf w}_{\tau_{j_2}}^{(i_2)}
\stackrel{\rm def}{=}J[\Phi]_{T,t}^{(2)},
\end{equation}

\vspace{2mm}
\noindent
where $\Phi(t_1,t_2):\ [t, T]^2\to\mathbb{R}$ is a nonrandom function 
(the properties of this function
will be specified further),
$\Delta{\bf w}_{\tau_{j}}^{(i)}=
{\bf w}_{\tau_{j+1}}^{(i)}-{\bf w}_{\tau_{j}}^{(i)}$
$(i=0, 1,\ldots,m)$,
$\left\{\tau_{j}\right\}_{j=0}^{N}$ is the partition 
of the interval $[t,T]$ satisfying the condition {\rm (\ref{1111})}.

Denote
\begin{equation}
\label{dom1}
D_2=\{(t_1,t_2):\ t\le t_1<t_2\le T\}.
\end{equation}

\vspace{3mm}

We will use the same symbol $D_2$ to denote the open and closed 
domains corresponding to the domain $D_2$ defined by (\ref{dom1}).
However, we always specify what domain we consider (open or closed).

Also we will write $\Phi(t_1,t_2)\in C(D_2)$
if 
$\Phi(t_1,t_2)$ is a continuous nonrandom function of two variables
in the closed domain $D_2$.

Let us consider the iterated Ito stochastic integral

\vspace{-1mm}
$$
I[\Phi]_{T,t}^{(2)}\stackrel{\rm def}{=}
\int\limits_t^T\int\limits_t^{t_2}
\Phi(t_1,t_2)d{\bf w}_{t_1}^{(i_1)}d{\bf w}_{t_2}^{(i_2)},
$$

\vspace{2mm}
\noindent
where $\Phi(t_1,t_2)\in C(D_2).$

\vspace{2mm}

{\bf Lemma 2}\ \cite{2013} (also see \cite{2017-1xx} and references therein).
{\it Suppose that
$\Phi(t_1,t_2)\in C(D_2)$ or $\Phi(t_1,t_2)$ 
is a continuous nonrandom function in the open domain $D_2$ and bounded at its boundary.
Then

\vspace{-1mm}
\begin{equation}
\label{30.52}
I[\Phi]_{T,t}^{(2)}=\hbox{\vtop{\offinterlineskip\halign{
\hfil#\hfil\cr
{\rm l.i.m.}\cr
$\stackrel{}{{}_{N\to \infty}}$\cr
}} }
\sum_{j_k=0}^{N-1}
\sum_{j_1=0}^{j_{2}-1}
\Phi(\tau_{j_1},\tau_{j_2})
\Delta {\bf w}_{\tau_{j_1}}^{(i_1)}\Delta {\bf w}_{\tau_{j_2}}^{(i_2)}\ \ \ 
\hbox{{\rm w. p. 1}},
\end{equation}

\vspace{3mm}
\noindent
where $\Delta{\bf w}_{\tau_{j}}^{(i)}=
{\bf w}_{\tau_{j+1}}^{(i)}-{\bf w}_{\tau_{j}}^{(i)}$
$(i=0, 1,\ldots,m)$,
$\left\{\tau_{j}\right\}_{j=0}^{N}$ is the partition 
of the interval $[t,T]$ satisfying the condition {\rm (\ref{1111})}.
}

\vspace{3mm}

{\bf Lemma 3}\ \cite{2013} (also see \cite{2017-1xx} and references therein).
{\it Suppose that every $\varphi_l(\tau)$
$(l=1,\ 2)$ is a continuous nonrandom function at the interval $[t, T]$.
Then

\begin{equation}
\label{30.39}
J[\varphi_1]_{T,t}J[\varphi_2]_{T,t}=J[\Phi]_{T,t}^{(2)}\ \ \ 
\hbox{w. p. {\rm 1}},
\end{equation}

\vspace{3mm}
\noindent
where 
$$
\Phi(t_1,t_2)=\varphi_1(t_1)\varphi_2(t_2),\ \ \ \
J[\varphi_l]_{T,t}
=\int\limits_t^T \varphi_l(\tau) d{\bf w}_{\tau}^{(i_l)}\ \ \ 
(l=1,\ 2)
$$

\vspace{2mm}
\noindent
and the stochastic integral $J[\Phi]_{T,t}^{(2)}$ 
is defined
by the equality
{\rm (\ref{30.34}),}\ $i_1, i_2=0, 1,\ldots,m.$
}

\vspace{2mm}

In accordance to the standard relations between
Stratonovich and Ito stochastic integrals we have w.~p.~1 \cite{KlPl2}

\vspace{-3mm}
\begin{equation}
\label{oop51}
J^{*}[\psi^{(2)}]_{T,t}=
J[\psi^{(2)}]_{T,t}+
\frac{1}{2}{\bf 1}_{\{i_1=i_2\ne 0\}}
\int\limits_t^T\psi_1(t_1)\psi_2(t_1)dt_1,
\end{equation}

\vspace{3mm}
\noindent
where ${\bf 1}_A$ is the indicator of the set $A$.

Let us define the function $K^{*}(t_1,t_2)$ at the square
$[t,T]^2$ as follows

\vspace{-1mm}
\begin{equation}
\label{1999.1}
K^{*}(t_1,t_2)=\psi_1(t_1)\psi_2(t_2)
\Biggl({\bf 1}_{\{t_1<t_{2}\}}+
\frac{1}{2}{\bf 1}_{\{t_1=t_{2}\}}\Biggr)
=K(t_1,t_2)+\frac{1}{2}{\bf 1}_{\{t_1=t_{2}\}}\psi_1(t_1)\psi_2(t_2),
\end{equation}

\vspace{2mm}
\noindent
where

\vspace{-2mm}
$$
K(t_1,t_2)=
\begin{cases}
\psi_1(t_1)\psi_2(t_2),\ &t_1<t_2\cr\cr
0,\ &\hbox{\rm otherwise}
\end{cases},\ \ \ t_1, t_2\in[t, T]
$$

\vspace{5mm}
\noindent
and ${\bf 1}_A$ is the indicator of the set $A$.

\vspace{2mm}

{\bf Lemma 4}\ 
\cite{2013} (also see \cite{2017-1xx} and references therein).
{\it Under the conditions of Theorem {\rm 1}
the following relation

\vspace{-1mm}
\begin{equation}
\label{30.36}
J[{K^{*}}]_{T,t}^{(2)}=
J^{*}[\psi^{(2)}]_{T,t}
\end{equation}

\vspace{2mm}
\noindent
is valid w.~p.~{\rm 1}, 
where $J[{K^{*}}]_{T,t}^{(2)}$ is defined by the equality 
{\rm (\ref{30.34})}.}

\vspace{2mm}

{\bf Proof.} Substituting
(\ref{1999.1}) into (\ref{30.34}) and using Lemmas 1 and 2,
it is easy to see that 

\vspace{-1mm}
\begin{equation}
\label{30.37}
J[{K^{*}}]_{T,t}^{(2)}
=
J[\psi^{(2)}]_{T,t}+
\frac{1}{2}{\bf 1}_{\{i_1=i_2\ne 0\}}
\int\limits_t^T\psi_1(t_1)\psi_2(t_1)dt_1
=J^{*}[\psi^{(2)}]_{T,t}\ \ \ \hbox{w.\ p.\ 1}.
\end{equation}

\vspace{3mm}

Let us consider the following generalized double Fourier sum

\vspace{1mm}
$$
\sum_{j_1=0}^{p_1}\sum_{j_2=0}^{p_2}
C_{j_2j_1}\phi_{j_1}(t_1)\phi_{j_2}(t_2),
$$

\vspace{5mm}
\noindent
where $C_{j_2j_1}$ is the Fourier coefficient of the form

\vspace{1mm}
\begin{equation}
\label{1}
C_{j_2j_1}=\int\limits_{[t,T]^2}
K^{*}(t_1,t_2)\phi_{j_1}(t_1)\phi_{j_2}(t_2)dt_1dt_2.
\end{equation}

\vspace{3mm}

Substitute the relation

\vspace{1mm}
$$
K^{*}(t_1,t_2)=
\sum_{j_1=0}^{p_1}\sum_{j_2=0}^{p_2}
C_{j_2j_1}\phi_{j_1}(t_1)\phi_{j_2}(t_2)+
K^{*}(t_1,t_2)
-\sum_{j_1=0}^{p_1}\sum_{j_2=0}^{p_2}
C_{j_2j_1}\phi_{j_1}(t_1)\phi_{j_2}(t_2)
$$

\vspace{4mm}
\noindent
with finite $p_1$ and $p_2$
into $J[{K^{*}}]_{T,t}^{(2)}.$ 
Then, using Lemma 3, we obtain

\vspace{1mm}
\begin{equation}
\label{proof1}
J^{*}[\psi^{(2)}]_{T,t}=
\sum_{j_1=0}^{p_1}\sum_{j_2=0}^{p_2}
C_{j_2j_1}
\zeta_{j_1}^{(i_1)}\zeta_{j_2}^{(i_2)}+
J[R_{p_1p_2}]_{T,t}^{(2)}\ \ \ \hbox{w. p. {\rm 1}},
\end{equation}

\vspace{4mm}
\noindent
where the stochastic integral
$J[R_{p_1p_2}]_{T,t}^{(2)}$
is defined in accordance with (\ref{30.34}) and

\vspace{1mm}
\begin{equation}
\label{30.46}
R_{p_1p_2}(t_1,t_2)=
K^{*}(t_1,t_2)-
\sum_{j_1=0}^{p_1}\sum_{j_2=0}^{p_2}
C_{j_2j_1}\phi_{j_1}(t_1)\phi_{j_2}(t_2),
\end{equation}

\vspace{1mm}
$$
\zeta_{j}^{(i)}=\int\limits_t^T \phi_{j}(\tau) d{\bf w}_{\tau}^{(i)},
$$

\vspace{3mm}

$$
J[R_{p_1p_2}]_{T,t}^{(2)}=\int\limits_t^T\int\limits_t^{t_2}
R_{p_1p_2}(t_1,t_2)d{\bf w}_{t_1}^{(i_1)}d{\bf w}_{t_2}^{(i_2)}
+\int\limits_t^T\int\limits_t^{t_1}
R_{p_1p_2}(t_1,t_2)d{\bf w}_{t_2}^{(i_2)}d{\bf w}_{t_1}^{(i_1)}+
$$

\vspace{1mm}
$$
+{\bf 1}_{\{i_1=i_2\ne 0\}}
\int\limits_t^T R_{p_1p_2}(t_1,t_1)dt_1.
$$

\vspace{5mm}

Let us consider the case $i_1, i_2\ne 0$ (another cases can be considered
absolutely analogously).
Using standard estimates for 
moments of stochastic integrals \cite{1}, we obtain 

\vspace{1mm}
$$
{\sf M}\left\{\left(J[R_{p_1p_2}]_{T,t}^{(2)}\right)^{2}
\right\}=
$$

\vspace{2mm}
$$
={\sf M}\left\{\left(\int\limits_t^T\int\limits_t^{t_2}
R_{p_1p_2}(t_1,t_2)d{\bf w}_{t_1}^{(i_1)}d{\bf w}_{t_2}^{(i_2)}
+\int\limits_t^T\int\limits_t^{t_1}
R_{p_1p_2}(t_1,t_2)d{\bf w}_{t_2}^{(i_2)}d{\bf w}_{t_1}^{(i_1)}
\right)^2\right\}+
$$

$$
+{\bf 1}_{\{i_1=i_2\ne 0\}}
\left(\int\limits_t^T R_{p_1p_2}(t_1,t_1)dt_1\right)^2\le
$$

\vspace{1mm}
$$
\le 
2\left(\int\limits_t^T\int\limits_t^{t_2}
\left(R_{p_1p_2}(t_1,t_2)\right)^{2}dt_1 dt_2
+
\int\limits_t^T\int\limits_t^{t_1}
\left(R_{p_1p_2}(t_1,t_2)\right)^{2}dt_2 dt_1\right)+
$$

$$
+
{\bf 1}_{\{i_1=i_2\ne 0\}}
\left(\int\limits_t^T R_{p_1p_2}(t_1,t_1)dt_1\right)^2=
$$

\begin{equation}
\label{newbegin1}
=
2\int\limits_{[t, T]^2}
\left(R_{p_1p_2}(t_1,t_2)\right)^{2}dt_1 dt_2
+
{\bf 1}_{\{i_1=i_2\ne 0\}}
\left(\int\limits_t^T R_{p_1p_2}(t_1,t_1)dt_1\right)^2.
\end{equation}

\vspace{5mm}

We have

\vspace{-1mm}
$$
\int\limits_{[t, T]^2}
\left(R_{p_1p_2}(t_1,t_2)\right)^{2}dt_1 dt_2=
$$

\vspace{1mm}
$$
=
\int\limits_{[t, T]^2}
\Biggl(
K^{*}(t_1,t_2)-
\sum_{j_1=0}^{p_1}\sum_{j_2=0}^{p_2}C_{j_2 j_1}
\phi_{j_1}(t_1)\phi_{j_2}(t_2)\Biggr)^2 dt_1 dt_2=
$$

\vspace{1mm}
$$
=\int\limits_{[t, T]^2}
\Biggl(
K(t_1,t_2)-
\sum_{j_1=0}^{p_1}\sum_{j_2=0}^{p_2}C_{j_2 j_1}
\phi_{j_1}(t_1)\phi_{j_2}(t_2)\Biggr)^2 dt_1 dt_2.
$$

\vspace{5mm}

The function $K(t_1,t_2)$ is piecewise continuous in the 
square $[t, T]^2$.
At this situation it is well-known that the generalized
multiple Fourier series 
of the function $K(t_1,t_2)\in L_2([t, T]^2)$ is converging 
to this function in the square $[t, T]^2$ in the mean-square sense, i.e.

$$
\hbox{\vtop{\offinterlineskip\halign{
\hfil#\hfil\cr
{\rm lim}\cr
$\stackrel{}{{}_{p_1,p_2\to \infty}}$\cr
}} }\Biggl\Vert
K(t_1,t_2)-
\sum_{j_1=0}^{p_1}\sum_{j_2=0}^{p_2}
C_{j_2 j_1}\prod_{l=1}^{2} \phi_{j_l}(t_l)\Biggr\Vert_{L_2([t,T]^2)}=0,
$$

\vspace{4mm}
\noindent
where
$$
\left\Vert f\right\Vert_{L_2([t,T]^2)}=\left(\int\limits_{[t,T]^2}
f^2(t_1,t_2)dt_1dt_2\right)^{1/2}.
$$

\vspace{3mm}

So, we obtain

\begin{equation}
\label{newbegin2}
\hbox{\vtop{\offinterlineskip\halign{
\hfil#\hfil\cr
{\rm lim}\cr
$\stackrel{}{{}_{p_1,p_2\to \infty}}$\cr
}} }
\int\limits_{[t, T]^2}
\left(R_{p_1p_2}(t_1,t_2)\right)^{2}dt_1 dt_2=0.
\end{equation}

\vspace{4mm}

Note that

\vspace{-2mm}
$$
\int\limits_t^T R_{p_1p_2}(t_1,t_1)dt_1=
$$

\vspace{1mm}
$$
=
\int\limits_t^T
\left(
\frac{1}{2}\psi_1(t_1)\psi_2(t_1) - 
\sum_{j_1=0}^{p_1}\sum_{j_2=0}^{p_2}C_{j_2 j_1}
\phi_{j_1}(t_1)\phi_{j_2}(t_1)\right) dt_1=
$$

\vspace{1mm}
$$
=
\frac{1}{2}\int\limits_t^T
\psi_1(t_1)\psi_2(t_1)dt_1 -
\sum_{j_1=0}^{p_1}\sum_{j_2=0}^{p_2}C_{j_2 j_1}
\int\limits_t^T\phi_{j_1}(t_1)\phi_{j_2}(t_1)dt_1=
$$

\vspace{1mm}
$$
=
\frac{1}{2}\int\limits_t^T
\psi_1(t_1)\psi_2(t_1)dt_1 -
\sum_{j_1=0}^{p_1}\sum_{j_2=0}^{p_2}C_{j_2 j_1}
{\bf 1}_{\{j_1=j_2\}}=
$$

\vspace{1mm}
\begin{equation}
\label{newbegin3}
=
\frac{1}{2}\int\limits_t^T
\psi_1(t_1)\psi_2(t_1)dt_1 -
\sum_{j_1=0}^{{\rm min}\{p_1,p_2\}}C_{j_1 j_1}.
\end{equation}

\vspace{5mm}

From (\ref{newbegin3}) we obtain

\vspace{-1mm}
\begin{equation}
\label{dds1}
\lim\limits_{p_1\to\infty}
\lim\limits_{p_2\to\infty}\int\limits_t^T R_{p_1p_2}(t_1,t_1)dt_1
=
\end{equation}

\vspace{1mm}
$$
=\frac{1}{2}\int\limits_t^T
\psi_1(t_1)\psi_2(t_1)dt_1 -
\lim\limits_{p_1\to\infty}
\sum_{j_1=0}^{p_1}C_{j_1 j_1}=
$$

\vspace{1mm}
$$
=\frac{1}{2}\int\limits_t^T
\psi_1(t_1)\psi_2(t_1)dt_1 -
\sum_{j_1=0}^{\infty}C_{j_1 j_1}
=
$$

\vspace{1mm}
\begin{equation}
\label{s1}
=\lim\limits_{p_1,p_2\to\infty}\int\limits_t^T R_{p_1p_2}(t_1,t_1)dt_1.
\end{equation}

\vspace{4mm}

Note that the existence of the limit 

\vspace{-1mm}
$$
\lim\limits_{p_1\to\infty}
\sum_{j_1=0}^{p_1}C_{j_1 j_1}
$$

\vspace{2mm}
\noindent
is proved in 
\cite{2017-1xx} (Sect.~2.1.1, 2.1.2) for the polynomial
and trigonometric cases.

If we prove the following relation

\vspace{-1mm}
\begin{equation}
\label{s11}
\lim\limits_{p_1\to\infty}
\lim\limits_{p_2\to\infty}\int\limits_t^T R_{p_1p_2}(t_1,t_1)dt_1=0,
\end{equation}

\vspace{2mm}
\noindent
then from (\ref{s1}) we get

\vspace{-1mm}
\begin{equation}
\label{44}
\frac{1}{2}\int\limits_t^T
\psi_1(t_1)\psi_2(t_1)dt_1 =
\sum_{j_1=0}^{\infty}C_{j_1 j_1},
\end{equation}

\begin{equation}
\label{444}
\lim\limits_{p_1,p_2\to\infty}\int\limits_t^T R_{p_1p_2}(t_1,t_1)dt_1=0.
\end{equation}
       
\vspace{4mm}

From (\ref{newbegin1}), (\ref{newbegin2}), and (\ref{444}) we obtain

\vspace{1mm}
$$
\lim\limits_{p_1,p_2\to\infty}
{\sf M}\left\{\left(J[R_{p_1p_2}]_{T,t}^{(2)}\right)^{2}
\right\}=0
$$

\vspace{4mm}
\noindent
and Theorem 1 will be proved.

Let us expand the function $K^{*}(t_1,t_2)$ 
(see (\ref{1999.1})) using the variable 
$t_1$, when $t_2$ is fixed, into the generalized Fourier series 
at the interval $(t, T)$

\begin{equation}
\label{leto8001yes1}
K^{*}(t_1,t_2)=
\sum_{j_1=0}^{\infty}C_{j_1}(t_2)\phi_{j_1}(t_1)\ \ \ (t_1\ne t, T),
\end{equation}

\vspace{4mm}
\noindent
where
$$
C_{j_1}(t_2)=\int\limits_t^T
K^{*}(t_1,t_2)\phi_{j_1}(t_1)dt_1=\psi_2(t_2)
\int\limits_t^{t_2}\psi_1(t_1)\phi_{j_1}(t_1)dt_1.
$$

\vspace{3mm}

The equality (\ref{leto8001yes1}) is 
satisfied
pointwise in each point of the interval $(t, T)$ with respect to the 
variable $t_1$, when $t_2\in [t, T]$ is fixed, due to a
piecewise
smoothness of the function $K^{*}(t_1,t_2)$ with respect to the variable 
$t_1\in [t, T]$ ($t_2$ is fixed). 

Note also that due to well-known properties of the Fourier--Legendre series 
and trigonometric Fourier series, 
the series (\ref{leto8001yes1}) converges when $t_1=t$ and $t_1=T$. 

Obtaining (\ref{leto8001yes1}), we also used the fact that the right-hand side
of (\ref{leto8001yes1}) converges when $t_1=t_2$ (point of a finite 
discontinuity
of the function $K(t_1,t_2)$) to the value

\vspace{1mm}
$$
\frac{1}{2}\left(K(t_2-0,t_2)+K(t_2+0,t_2)\right)=
\frac{1}{2}\psi_1(t_2)\psi_2(t_2)=
K^{*}(t_2,t_2).
$$

\vspace{4mm}

The function $C_{j_1}(t_2)$ is a continuously differentiable
one at the interval $[t, T]$. 
Let us expand it into the generalized Fourier series 
at the interval $(t, T)$

\begin{equation}
\label{leto8002yes}
C_{j_1}(t_2)=
\sum_{j_2=0}^{\infty}C_{j_2 j_1}\phi_{j_2}(t_2)\ \ \ (t_2\ne t, T),
\end{equation}

\vspace{4mm}
\noindent
where
$$
C_{j_2 j_1}=\int\limits_t^T
C_{j_1}(t_2)\phi_{j_2}(t_2)dt_2=
\int\limits_t^T
\psi_2(t_2)\phi_{j_2}(t_2)\int\limits_t^{t_2}
\psi_1(t_1)\phi_{j_1}(t_1)dt_1 dt_2
$$

\vspace{2mm}
\noindent
and the equality (\ref{leto8002yes}) is satisfied pointwise at any point 
of the interval $(t, T)$ (the right-hand side 
of (\ref{leto8002yes}) converges 
when $t_2=t$ and $t_1=T$).

Let us substitute (\ref{leto8002yes}) into (\ref{leto8001yes1})

\begin{equation}
\label{leto8003yes}
K^{*}(t_1,t_2)=
\sum_{j_1=0}^{\infty}\sum_{j_2=0}^{\infty}C_{j_2 j_1}
\phi_{j_1}(t_1)\phi_{j_2}(t_2),\ \ \ (t_1, t_2)\in (t, T)^2.
\end{equation}

\vspace{4mm}

Futhermore,
the series on the right-hand side of (\ref{leto8003yes}) converges at the 
boundary
of the square $[t, T]^2$.

From (\ref{30.46}) and (\ref{leto8003yes}) we obtain

\begin{equation}
\label{y1}
\lim\limits_{p_1\to\infty}
\lim\limits_{p_2\to\infty}R_{p_1p_2}(t_1,t_1)=0\ \ \ \hbox{when}\ \ \ t_1\in (t, T).
\end{equation}

\vspace{4mm}

Since the integral
\begin{equation}
\label{1113}
\int\limits_t^T R_{p_1p_2}(t_1,t_1)dt_1
\end{equation}

\vspace{2mm}
\noindent
exists as Riemann integral, then this integral equals to the 
corresponding Lebesgue integral. 
Moreover, 

\begin{equation}
\label{1110}
\lim\limits_{p_1\to\infty}\lim\limits_{p_2\to\infty}
R_{p_1p_2}(t_1,t_1)=0\ \ \ \hbox{when}\ \ \
t_1\in (t, T),
\end{equation}

\vspace{4mm}
\noindent
where the left-hand side of (\ref{1110}) is bounded on $[t, T].$

According to (\ref{30.46}), (\ref{leto8001yes1})--(\ref{leto8003yes}), we have

$$
R_{p_1p_2}(t_1,t_2)=\left(K^{*}(t_1,t_2)-\sum\limits_{j_1=0}^{p_1}
C_{j_1}(t_2)\phi_{j_1}(t_1)\right)+
$$

\vspace{2mm}
$$
+\left(
\sum\limits_{j_1=0}^{p_1}\left(C_{j_1}(t_2)-
\sum\limits_{j_2=0}^{p_2}
C_{j_2j_1}\phi_{j_2}(t_2)\right)
\phi_{j_1}(t_1)\right).
$$

\vspace{5mm}

Then, appling two times (we mean here an iterated passage to the limit
$\lim\limits_{p_1\to\infty}\lim\limits_{p_2\to\infty}$)
the Lebesgue's Dominated Convergence Theorem to the integral
(\ref{1113}), 
we obtain 

\vspace{1mm}
$$
\lim\limits_{p_1\to\infty}\lim\limits_{p_2\to\infty}
\int\limits_t^T
R_{p_1p_2}(t_1,t_1)dt_1=0.
$$

\vspace{3mm}

For a discussion of the choice of integrable majorants
when applying Lebesgue's 
Dominated Convergence Theorem to the integral (\ref{1113})
for the polynomial and trigonometric cases, see \cite{2017-1xx}
(Sect.~2.4.1), \cite{xxxxx} (Sect.~2).

Note that the developement of the approach from this article 
can be found in 
\cite{2017-1xx} (Sect.~2.4), \cite{xxxxx}.

\vspace{5mm}

\section{Some Recent Results on Expansion
of Iterated Stratonovich Stochastic Integrals
of Multiplicities 2 to 6}

\vspace{5mm}

Recently, a new approach to the expansion and mean-square 
approximation of iterated Stratonovich stochastic integrals has been obtained
\cite{2017-1xx} (Sect.~2.10--2.16), \cite{32} (Sect.~13--19), 
\cite{15a} (Sect.~5--11), \cite{arxiv-11} (Sect.~7--13), \cite{new-art-1-xxy}
(Sect.~4--9).
Let us formulate four theorems that were obtained using this approach.

\vspace{2mm}

{\bf Theorem 2}\ \cite{2017-1xx}, \cite{32}, \cite{15a}, \cite{arxiv-11}, \cite{new-art-1-xxy}.\
{\it Suppose 
that $\{\phi_j(x)\}_{j=0}^{\infty}$ is a complete orthonormal system of 
Legendre polynomials or trigonometric functions in the space $L_2([t, T]).$
Furthermore, let $\psi_1(\tau), \psi_2(\tau),$ $\psi_3(\tau)$ are continuously dif\-ferentiable 
nonrandom functions on $[t, T].$ 
Then, for the 
iterated Stra\-to\-no\-vich stochastic integral of third multiplicity

$$
J^{*}[\psi^{(3)}]_{T,t}={\int\limits_t^{*}}^T\psi_3(t_3)
{\int\limits_t^{*}}^{t_3}\psi_2(t_2)
{\int\limits_t^{*}}^{t_2}\psi_1(t_1)
d{\bf w}_{t_1}^{(i_1)}
d{\bf w}_{t_2}^{(i_2)}d{\bf w}_{t_3}^{(i_3)}\ \ \ (i_1,i_2,i_3=0,1,\ldots,m)
$$

\vspace{4mm}
\noindent
the following 
relations

\vspace{-1mm}
\begin{equation}
\label{fin1}
J^{*}[\psi^{(3)}]_{T,t}
=\hbox{\vtop{\offinterlineskip\halign{
\hfil#\hfil\cr
{\rm l.i.m.}\cr
$\stackrel{}{{}_{p\to \infty}}$\cr
}} }
\sum\limits_{j_1, j_2, j_3=0}^{p}
C_{j_3 j_2 j_1}\zeta_{j_1}^{(i_1)}\zeta_{j_2}^{(i_2)}\zeta_{j_3}^{(i_3)},
\end{equation}

\vspace{3mm}
\begin{equation}
\label{fin2}
{\sf M}\left\{\left(
J^{*}[\psi^{(3)}]_{T,t}-
\sum\limits_{j_1, j_2, j_3=0}^{p}
C_{j_3 j_2 j_1}\zeta_{j_1}^{(i_1)}\zeta_{j_2}^{(i_2)}\zeta_{j_3}^{(i_3)}\right)^2\right\}
\le \frac{C}{p}
\end{equation}

\vspace{5mm}
\noindent
are fulfilled, where $i_1, i_2, i_3=0,1,\ldots,m$ in {\rm (\ref{fin1})} and 
$i_1, i_2, i_3=1,\ldots,m$ in {\rm (\ref{fin2})},
constant $C$ is independent of $p,$

$$
C_{j_3 j_2 j_1}=\int\limits_t^T\psi_3(t_3)\phi_{j_3}(t_3)
\int\limits_t^{t_3}\psi_2(t_2)\phi_{j_2}(t_2)
\int\limits_t^{t_2}\psi_1(t_1)\phi_{j_1}(t_1)dt_1dt_2dt_3
$$

\vspace{4mm}
\noindent
and
$$
\zeta_{j}^{(i)}=
\int\limits_t^T \phi_{j}(\tau) d{\bf f}_{\tau}^{(i)}
$$ 

\vspace{2mm}
\noindent
are independent standard Gaussian random variables for various 
$i$ or $j$ {\rm (}in the case when $i\ne 0${\rm );}
another notations are the same as in Theorem~{\rm 1}.}

\vspace{2mm}

{\bf Theorem 3}\ \cite{2017-1xx}, \cite{32}, \cite{15a}, \cite{arxiv-11}, \cite{new-art-1-xxy}.\ 
{\it Let
$\{\phi_j(x)\}_{j=0}^{\infty}$ be a complete orthonormal system of 
Legendre polynomials or trigonometric functions in the space $L_2([t, T]).$
Furthermore, let $\psi_1(\tau), \ldots,$ $\psi_4(\tau)$ be continuously dif\-ferentiable 
nonrandom functions on $[t, T].$ 
Then, for the 
iterated Stra\-to\-no\-vich stochastic integral of fourth multiplicity

\begin{equation}
\label{fin0}
J^{*}[\psi^{(4)}]_{T,t}={\int\limits_t^{*}}^T\psi_4(t_4)
{\int\limits_t^{*}}^{t_4}\psi_3(t_3)
{\int\limits_t^{*}}^{t_3}\psi_2(t_2)
{\int\limits_t^{*}}^{t_2}\psi_1(t_1)
d{\bf w}_{t_1}^{(i_1)}
d{\bf w}_{t_2}^{(i_2)}d{\bf w}_{t_3}^{(i_3)}d{\bf w}_{t_4}^{(i_4)}
\end{equation}

\vspace{4mm}
\noindent
the following 
relations

\begin{equation}
\label{fin3}
J^{*}[\psi^{(4)}]_{T,t}
=\hbox{\vtop{\offinterlineskip\halign{
\hfil#\hfil\cr
{\rm l.i.m.}\cr
$\stackrel{}{{}_{p\to \infty}}$\cr
}} }
\sum\limits_{j_1, j_2, j_3,j_4=0}^{p}
C_{j_4j_3 j_2 j_1}\zeta_{j_1}^{(i_1)}\zeta_{j_2}^{(i_2)}\zeta_{j_3}^{(i_3)}\zeta_{j_4}^{(i_4)},
\end{equation}

\vspace{3mm}

\begin{equation}
\label{fin4}
{\sf M}\left\{\left(
J^{*}[\psi^{(4)}]_{T,t}-
\sum\limits_{j_1, j_2, j_3, j_4=0}^{p}
C_{j_4 j_3 j_2 j_1}\zeta_{j_1}^{(i_1)}\zeta_{j_2}^{(i_2)}\zeta_{j_3}^{(i_3)}
\zeta_{j_4}^{(i_4)}
\right)^2\right\}
\le \frac{C}{p^{1-\varepsilon}}
\end{equation}

\vspace{5mm}
\noindent
are fulfilled, where $i_1, \ldots , i_4=0,1,\ldots,m$ in {\rm (\ref{fin0}),} {\rm (\ref{fin3})} 
and $i_1, \ldots, i_4=1,\ldots,m$ in {\rm (\ref{fin4}),}
constant $C$ does not depend on $p,$
$\varepsilon$ is an arbitrary
small positive real number 
for the case of complete orthonormal system of 
Legendre polynomials in the space $L_2([t, T])$
and $\varepsilon=0$ for the case of
complete orthonormal system of 
trigonometric functions in the space $L_2([t, T]),$

$$
C_{j_4 j_3 j_2 j_1}=
$$

$$
=
\int\limits_t^T\psi_4(t_4)\phi_{j_4}(t_4)
\int\limits_t^{t_4}\psi_3(t_3)\phi_{j_3}(t_3)
\int\limits_t^{t_3}\psi_2(t_2)\phi_{j_2}(t_2)
\int\limits_t^{t_2}\psi_1(t_1)\phi_{j_1}(t_1)dt_1dt_2dt_3dt_4;
$$

\vspace{4mm}
\noindent
another notations are the same as in Theorem~{\rm 2}.}

\vspace{2mm}

{\bf Theorem 4}\ \cite{2017-1xx}, \cite{32}, \cite{15a}, \cite{arxiv-11}, \cite{new-art-1-xxy}.\
{\it Assume 
that $\{\phi_j(x)\}_{j=0}^{\infty}$ is a complete orthonormal system of 
Legendre polynomials or trigonometric functions in the space $L_2([t, T])$
and $\psi_1(\tau), \ldots,$ $\psi_5(\tau)$ are continuously dif\-ferentiable 
nonrandom functions on $[t, T].$ 
Then, for the 
iterated Stra\-to\-no\-vich stochastic integral of fifth multiplicity

\begin{equation}
\label{fin7}
J^{*}[\psi^{(5)}]_{T,t}={\int\limits_t^{*}}^T\psi_5(t_5)
\ldots
{\int\limits_t^{*}}^{t_2}\psi_1(t_1)
d{\bf w}_{t_1}^{(i_1)}
\ldots d{\bf w}_{t_5}^{(i_5)}
\end{equation}

\vspace{4mm}
\noindent
the following 
relations

\begin{equation}
\label{fin8}
J^{*}[\psi^{(5)}]_{T,t}
=\hbox{\vtop{\offinterlineskip\halign{
\hfil#\hfil\cr
{\rm l.i.m.}\cr
$\stackrel{}{{}_{p\to \infty}}$\cr
}} }
\sum\limits_{j_1,\ldots,j_5=0}^{p}
C_{j_5 \ldots j_1}\zeta_{j_1}^{(i_1)}\ldots \zeta_{j_5}^{(i_5)},
\end{equation}

\vspace{3mm}

\begin{equation}
\label{fin9}
{\sf M}\left\{\left(
J^{*}[\psi^{(5)}]_{T,t}-
\sum\limits_{j_1, \ldots, j_5=0}^{p}
C_{j_5 \ldots j_1}\zeta_{j_1}^{(i_1)}\ldots
\zeta_{j_5}^{(i_5)}
\right)^2\right\}
\le \frac{C}{p^{1-\varepsilon}}
\end{equation}

\vspace{5mm}
\noindent
are fulfilled, where $i_1, \ldots , i_5=0,1,\ldots,m$ in {\rm (\ref{fin7}),} {\rm (\ref{fin8})} 
and $i_1, \ldots, i_5=1,\ldots,m$ in {\rm (\ref{fin9}),}
constant $C$ is independent of $p,$
$\varepsilon$ is an arbitrary
small positive real number 
for the case of complete orthonormal system of 
Legendre polynomials in the space $L_2([t, T])$
and $\varepsilon=0$ for the case of
complete orthonormal system of 
trigonometric functions in the space $L_2([t, T]),$

$$
C_{j_5 \ldots j_1}=
\int\limits_t^T\psi_5(t_5)\phi_{j_5}(t_5)\ldots
\int\limits_t^{t_2}\psi_1(t_1)\phi_{j_1}(t_1)dt_1\ldots dt_5;
$$

\vspace{3mm}
\noindent
another notations are the same as in Theorems~{\rm 2, 3}.}

\vspace{2mm}

{\bf Theorem 5}\ \cite{2017-1xx}, \cite{32}, \cite{15a}, \cite{arxiv-11}.\
{\it Suppose that 
$\{\phi_j(x)\}_{j=0}^{\infty}$ is a complete orthonormal system of 
Legendre polynomials or trigonometric functions in the space $L_2([t, T]).$
Then, for the 
iterated Stra\-to\-no\-vich stochastic integral of sixth multiplicity

\begin{equation}
\label{after10001qu1}
J_{T,t}^{*(i_1\ldots i_6)}={\int\limits_t^{*}}^T
\ldots
{\int\limits_t^{*}}^{t_2}
d{\bf w}_{t_1}^{(i_1)}
\ldots d{\bf w}_{t_6}^{(i_6)}
\end{equation}

\vspace{3mm}
\noindent
the following 
expansion 

\vspace{-1mm}
$$
J_{T,t}^{*(i_1\ldots i_6)}
=\hbox{\vtop{\offinterlineskip\halign{
\hfil#\hfil\cr
{\rm l.i.m.}\cr
$\stackrel{}{{}_{p\to \infty}}$\cr
}} }
\sum\limits_{j_1, \ldots, j_6=0}^{p}
C_{j_6 \ldots j_1}\zeta_{j_1}^{(i_1)}\ldots
\zeta_{j_6}^{(i_6)}
$$

\vspace{4mm}
\noindent
that converges in the mean-square sense is valid, where
$i_1, \ldots, i_6=0, 1,\ldots,m,$

$$
C_{j_6 \ldots j_1}=
\int\limits_t^T\phi_{j_6}(t_6)\ldots
\int\limits_t^{t_2}\phi_{j_1}(t_1)dt_1\ldots dt_6;
$$

\vspace{3mm}
\noindent
another notations are the same as in Theorems~{\rm 2--4}.}

\vspace{2mm}

Recently the equality (\ref{44}) was proved in \cite{rybakov7000x} (also see \cite{2017-1xx} (Sect.~2.1.4))
for the case of an arbitrary complete orthonormal system of 
functions in $L_2([t, T])$ and $\psi_1(\tau), \psi_2(\tau)\in L_2([t, T]).$
This means that we have the following generalizaion of Theorem~1.

\vspace{2mm}

{\bf Theorem~6}\ \cite{2017-1xx} (Sect.~2.1.4).\
{\it Suppose that 
$\{\phi_j(x)\}_{j=0}^{\infty}$ is an arbitrary complete orthonormal system of 
functions in the space $L_2([t, T])$
and $\psi_1(\tau), \psi_2(\tau)$ are continuous
functions on $[t, T]$. 
Then the iterated Stratonovich stochastic integral of the second multiplicity

\vspace{-1mm}
$$
J^{*}[\psi^{(2)}]_{T,t}={\int\limits_t^{*}}^T\psi_2(t_2)
{\int\limits_t^{*}}^{t_2}\psi_1(t_1)d{\bf w}_{t_1}^{(i_1)}
d{\bf w}_{t_2}^{(i_2)}\ \ \ (i_1, i_2=0, 1,\ldots,m)
$$

\vspace{2mm}
\noindent
is expanded into the 
multiple series

\vspace{-1mm}
$$
J^{*}[\psi^{(2)}]_{T,t}=
\hbox{\vtop{\offinterlineskip\halign{
\hfil#\hfil\cr
{\rm l.i.m.}\cr
$\stackrel{}{{}_{p_1, p_2\to \infty}}$\cr
}} }
\sum\limits_{j_1=0}^{p_1}
\sum\limits_{j_2=0}^{p_2}
C_{j_2j_1}\zeta_{j_1}^{(i_1)}\zeta_{j_2}^{(i_2)}
$$

\vspace{3mm}
\noindent
that converges 
in the mean-square sense{\rm ;} where
notations are the same as in Theorem~{\rm 1.}}

\vspace{2mm}

The condition of continuity of the functions
$\psi_1(\tau), \psi_2(\tau)$ 
is related to the definition 
of the Stratonovich stochastic integral that we use
(see \cite{KlPl2}).

\vspace{5mm}

\section{Theorems 1--6 from Point
of View of the Wong--Zakai Approximation}

\vspace{5mm}

The iterated Ito stochastic integrals and solutions
of Ito SDEs are complex and important functionals
from the independent components ${\bf f}_{s}^{(i)},$
$i=1,\ldots,m$ of the multidimensional
Wiener process ${\bf f}_{s},$ $s\in[0, T].$
Let ${\bf f}_{s}^{(i)p},$ $p\in\mathbb{N}$ 
be some approximation of
${\bf f}_{s}^{(i)},$
$i=1,\ldots,m$.
Suppose that 
${\bf f}_{s}^{(i)p}$
converges to
${\bf f}_{s}^{(i)},$
$i=1,\ldots,m$ if $p\to\infty$ in some sense and has
differentiable sample trajectories.

A natural question arises: if we replace 
${\bf f}_{s}^{(i)}$
by ${\bf f}_{s}^{(i)p},$
$i=1,\ldots,m$ in the functionals
mentioned above, will the resulting
functionals converge to the original
functionals from the components 
${\bf f}_{s}^{(i)},$
$i=1,\ldots,m$ of the multidimentional
Wiener process ${\bf f}_{s}$?
The answere to this question is negative 
in the general case. However, 
in the pioneering works of Wong E. and Zakai M. \cite{W-Z-1},
\cite{W-Z-2},
it was shown that under the special conditions and 
for some types of approximations 
of the Wiener process the answere is affirmative
with one peculiarity: the convergence takes place 
to the iterated Stratonovich stochastic integrals
and solutions of Stratonovich SDEs and not to iterated 
Ito stochastic integrals and solutions
of Ito SDEs.
The piecewise 
linear approximation 
as well as the regularization by convolution 
\cite{W-Z-1}-\cite{Watanabe} relate the 
mentioned types of approximations
of the Wiener process. The above approximation 
of stochastic integrals and solutions of SDEs 
is often called the Wong--Zakai approximation.

Let ${\bf w}_{\tau},$ $\tau\in[0, T]$ is a random vector with 
an $m+1$ components: ${\bf w}_{\tau}^{(i)}={\bf f}_{\tau}^{(i)}$ 
for $i=1,\ldots,m$ and 
${\bf w}_{\tau}^{(0)}=\tau,$\ 
${\bf f}_{\tau}^{(i)}$ $(i=1,\ldots,m)$
are independent standard Wiener processes.

It is well known that the following representation 
takes place \cite{Lipt}, \cite{7e}

\begin{equation}
\label{um1x}
{\bf w}_{\tau}^{(i)}-{\bf w}_{t}^{(i)}=
\sum_{j=0}^{\infty}\int\limits_t^{\tau}
\phi_j(s)ds\ \zeta_j^{(i)},\ \ \ \zeta_j^{(i)}=
\int\limits_t^T \phi_j(\tau)d{\bf w}_{\tau}^{(i)},
\end{equation}

\vspace{4mm}
\noindent
where $\tau\in[t, T],$ $t\ge 0,$
$\{\phi_j(x)\}_{j=0}^{\infty}$ is an arbitrary complete 
orthonormal system of functions in the space $L_2([t, T]),$ and
$\zeta_j^{(i)}$ are independent standard Gaussian 
random variables for various $i$ or $j.$
Moreover, the series (\ref{um1x}) converges for any $\tau\in [t, T]$
in the mean-square sense.

Let ${\bf w}_{\tau}^{(i)p}-{\bf w}_{t}^{(i)p}$ be 
the mean-square approximation of the process
${\bf w}_{\tau}^{(i)}-{\bf w}_{t}^{(i)},$
which has the following form

\vspace{-3mm}
\begin{equation}
\label{um1xx}
{\bf w}_{\tau}^{(i)p}-{\bf w}_{t}^{(i)p}=
\sum_{j=0}^{p}\int\limits_t^{\tau}
\phi_j(s)ds\ \zeta_j^{(i)}.
\end{equation}

\vspace{3mm}

From (\ref{um1xx}) we obtain

\vspace{-4mm}
\begin{equation}
\label{um1xxx}
d{\bf w}_{\tau}^{(i)p}=
\sum_{j=0}^{p}
\phi_j(\tau)\zeta_j^{(i)} d\tau.
\end{equation}

\vspace{4mm}

Consider the following iterated Riemann--Stieltjes
integral

\begin{equation}
\label{um1xxxx}
\int\limits_t^T
\psi_k(t_k)\ldots \int\limits_t^{t_2}\psi_1(t_1)
d{\bf w}_{t_1}^{(i_1)p_1}\ldots d{\bf w}_{t_k}^{(i_k)p_k},
\end{equation}

\vspace{4mm}
\noindent
where $i_1,\ldots,i_k=0,1,\ldots,m,$\ \ $p_1,\ldots,p_k\in\mathbb{N},$ 

\begin{equation}
\label{um1xxx1}
d{\bf w}_{\tau}^{(i)p}=
\left\{\begin{matrix}
d{\bf f}_{\tau}^{(i)p}\ &\hbox{\rm for}\ \ \ i=1,\ldots,m\cr\cr\cr
d\tau^p\ &\hbox{\rm for}\ \ \ i=0
\end{matrix}
,\right.
\end{equation}

\vspace{4mm}
\noindent
and $d{\bf f}_{\tau}^{(i)p},$ $d\tau^p$ are defined by the relation (\ref{um1xxx}).

Let us substitute (\ref{um1xxx}) into (\ref{um1xxxx})

\begin{equation}
\label{um1xxxx1}
\int\limits_t^T
\psi_k(t_k)\ldots \int\limits_t^{t_2}\psi_1(t_1)
d{\bf w}_{t_1}^{(i_1)p_1}\ldots d{\bf w}_{t_k}^{(i_k)p_k}=
\sum\limits_{j_1=0}^{p_1} \ldots \sum\limits_{j_k=0}^{p_k}
C_{j_k \ldots j_1}\prod\limits_{l=1}^k \zeta_{j_l}^{(i_l)},
\end{equation}

\vspace{4mm}
\noindent
where 
$$
\zeta_j^{(i)}=\int\limits_t^T \phi_j(\tau)d{\bf w}_{\tau}^{(i)}
$$ 

\vspace{2mm}
\noindent
are independent standard Gaussian random variables for various 
$i$ or $j$ (in the case when $i\ne 0$),
${\bf w}_{s}^{(i)}={\bf f}_{s}^{(i)}$ for
$i=1,\ldots,m$ and 
${\bf w}_{s}^{(0)}=s,$

$$
C_{j_k \ldots j_1}=\int\limits_t^T\psi_k(t_k)\phi_{j_k}(t_k)\ldots
\int\limits_t^{t_2}
\psi_1(t_1)\phi_{j_1}(t_1)
dt_1\ldots dt_k
$$

\vspace{4mm}
\noindent
is the Fourier coefficient.

Consider the following iterated
Stratonovich stochastic integrals

\vspace{-1mm}
\begin{equation}
\label{str1}
J^{*}[\psi^{(k)}]_{T,t}=
{\int\limits_t^{*}}^T\psi_k(t_k)\ldots {\int\limits_t^{*}}^{t_2}
\psi_1(t_1) d{\bf w}_{t_1}^{(i_1)}\ldots d{\bf w}_{t_k}^{(i_k)},
\end{equation}

\vspace{2mm}
\noindent
where every $\psi_l(\tau)$ $(l=1,\ldots,k)$ is a continuously differentiable nonrandom function 
at the interval $[t,T];$ another notations are the same as in (\ref{str}).

To best of our knowledge \cite{W-Z-1}-\cite{Watanabe}
the approximations of the Wiener process
in the Wong--Zakai approximation must satisfy fairly strong
restrictions
\cite{Watanabe}
(see Definition 7.1, pp.~480--481).
Moreover, approximations of the Wiener process that are
similar to (\ref{um1xx})
were not considered in \cite{W-Z-1}, \cite{W-Z-2}
(also see \cite{Watanabe}, Theorems 7.1, 7.2).
Therefore, the proof of analogs of Theorems 7.1 and 7.2 \cite{Watanabe}
for approximations of the Wiener 
process based on its series expansion (\ref{um1x})
should be carried out separately.
Thus, the mean-square convergence of the right-hand side
of (\ref{um1xxxx1}) to the appropriate 
iterated Stratonovich stochastic integral (\ref{str1})
does not follow from the results of the papers
\cite{W-Z-1}, \cite{W-Z-2} (also see \cite{Watanabe},
Theorems 7.1, 7.2).

However, in 
\cite{KlPl2}
(Sect.~5.8, pp.~202--204), \cite{KPS} (pp.~82-84),
\cite{KPW} (pp.~438-439),  
\cite{Zapad-9} (pp.~263-264) the authors use 
(without rigorous proof)
the Wong--Zakai approximation 
\cite{W-Z-1}-\cite{Watanabe} within the frames of the 
method of approximation of iterated Stratonovich
stochastic integrals 
based on the Karhunen--Loeve expansion of the Brownian bridge
process \cite{Mi2}.

From the other hand, Theorems 1--6 from this paper
can be considered as the proof of the
Wong--Zakai approximation for the iterated 
Stratonovich stochastic integrals (\ref{str1}) of multiplicities 2 to 6
based on the approximation (\ref{um1xx}) of the Wiener process.
At that, the Riemann--Stieltjes integrals (\ref{um1xxxx})
of multiplicities 2 to 6
converge in the mean-square sense 
to the appropriate Stratonovich 
stochastic integrals (\ref{str1}). 
Recall that
$\{\phi_j(x)\}_{j=0}^{\infty}$ (see (\ref{um1x}), (\ref{um1xx}), and
Theorems 1--5)
is a complete 
orthonormal system of Legendre polynomials or 
trigonometric functions 
in the space $L_2([t, T])$.

To illustrate the above reasoning, 
consider two examples for the case $k=2,$
$\psi_1(\tau),$ $\psi_2(\tau)\equiv 1;$ $i_1, i_2=1,\ldots,m.$

The first example relates to the piecewise linear approximation
of the multidimensional Wiener process (these approximations 
were considered in \cite{W-Z-1}-\cite{Watanabe}).

Let ${\bf b}_{\Delta}^{(i)}(t),$ $t\in[0, T]$ be the piecewise
linear approximation of the $i$th component ${\bf f}_t^{(i)}$
of the multidimensional standard Wiener process ${\bf f}_t,$
$t\in [0, T]$ with independent components
${\bf f}_t^{(i)},$ $i=1,\ldots,m,$ i.e.

\vspace{-2mm}
$$
{\bf b}_{\Delta}^{(i)}(t)={\bf f}_{k\Delta}^{(i)}+
\frac{t-k\Delta}{\Delta}\Delta{\bf f}_{k\Delta}^{(i)},
$$

\vspace{3mm}
\noindent
where 

\vspace{-2mm}
$$
\Delta{\bf f}_{k\Delta}^{(i)}={\bf f}_{(k+1)\Delta}^{(i)}-
{\bf f}_{k\Delta}^{(i)},\ \ \
t\in[k\Delta, (k+1)\Delta),\ \ \ k=0, 1,\ldots, N-1.
$$

\vspace{4mm}

Note that w.~p.~1

\vspace{-1mm}
\begin{equation}
\label{pridum}
\frac{d{\bf b}_{\Delta}^{(i)}}{dt}(t)=
\frac{\Delta{\bf f}_{k\Delta}^{(i)}}{\Delta},\ \ \
t\in[k\Delta, (k+1)\Delta),\ \ \ k=0, 1,\ldots, N-1.
\end{equation}

\vspace{4mm}

Consider the following iterated Riemann--Stieltjes
integral

\vspace{1mm}
$$
\int\limits_0^T
\int\limits_0^{s}
d{\bf b}_{\Delta}^{(i_1)}(\tau)d{\bf b}_{\Delta}^{(i_2)}(s),\ \ \ 
i_1,i_2=1,\ldots,m.
$$

\vspace{4mm}

Using (\ref{pridum}) and additive property of the Riemann--Stieltjes integral, 
we can write w.~p.~1

\vspace{2mm}
$$
\int\limits_0^T
\int\limits_0^{s}
d{\bf b}_{\Delta}^{(i_1)}(\tau)d{\bf b}_{\Delta}^{(i_2)}(s)=
\int\limits_0^T
\int\limits_0^{s}
\frac{d{\bf b}_{\Delta}^{(i_1)}}{d\tau}(\tau)d\tau
\frac{d {\bf b}_{\Delta}^{(i_2)}}{d s}(s)
ds =
$$

\vspace{3mm}
$$
=
\sum\limits_{l=0}^{N-1}\int\limits_{l\Delta}^{(l+1)\Delta}
\left(
\sum\limits_{q=0}^{l-1}\int\limits_{q\Delta}^{(q+1)\Delta}
\frac{\Delta{\bf f}_{q\Delta}^{(i_1)}}{\Delta}d\tau+
\int\limits_{l\Delta}^{s}
\frac{\Delta{\bf f}_{l\Delta}^{(i_1)}}{\Delta}d\tau\right)
\frac{\Delta{\bf f}_{l\Delta}^{(i_2)}}{\Delta}ds=
$$

\vspace{3mm}
$$
=\sum\limits_{l=0}^{N-1}\sum\limits_{q=0}^{l-1}
\Delta{\bf f}_{q\Delta}^{(i_1)}
\Delta{\bf f}_{l\Delta}^{(i_2)}+
\frac{1}{\Delta^2}\sum\limits_{l=0}^{N-1}
\Delta{\bf f}_{l\Delta}^{(i_1)}
\Delta{\bf f}_{l\Delta}^{(i_2)}
\int\limits_{l\Delta}^{(l+1)\Delta}
\int\limits_{l\Delta}^{s}d\tau ds=
$$

\vspace{3mm}
\begin{equation}
\label{oh-ty}
=\sum\limits_{l=0}^{N-1}\sum\limits_{q=0}^{l-1}
\Delta{\bf f}_{q\Delta}^{(i_1)}
\Delta{\bf f}_{l\Delta}^{(i_2)}+
\frac{1}{2}\sum\limits_{l=0}^{N-1}
\Delta{\bf f}_{l\Delta}^{(i_1)}
\Delta{\bf f}_{l\Delta}^{(i_2)}.
\end{equation}

\vspace{6mm}

Using (\ref{oh-ty}), it 
is not difficult to show 
that

\vspace{1mm}
\begin{equation}
\label{uh-111}
\hbox{\vtop{\offinterlineskip\halign{
\hfil#\hfil\cr
{\rm l.i.m.}\cr
$\stackrel{}{{}_{N\to \infty}}$\cr
}} }
\int\limits_0^T
\int\limits_0^{s}
d{\bf b}_{\Delta}^{(i_1)}(\tau)d{\bf b}_{\Delta}^{(i_2)}(s)=
\int\limits_0^T
\int\limits_0^{s}
d{\bf f}_{\tau}^{(i_1)}d{\bf f}_{s}^{(i_2)}+
\frac{1}{2}{\bf 1}_{\{i_1=i_2\}}\int\limits_0^T ds=
{\int\limits_0^{*}}^T
{\int\limits_0^{*}}^s
d{\bf f}_{\tau}^{(i_1)}d{\bf f}_{s}^{(i_2)},
\end{equation}

\vspace{5mm}
\noindent
where $\Delta\to 0$ if $N\to\infty$ ($N\Delta=T$).

Obviously, (\ref{uh-111}) agrees with Theorem 7.1 (see \cite{Watanabe},
p.~486).

The next example relates to the approximation
of the Wiener process based on its series expansion
(\ref{um1x}) for $t=0$, where
$\{\phi_j(x)\}_{j=0}^{\infty}$ 
is a complete 
orthonormal system of Legendre polynomials or 
trigonometric functions 
in the space $L_2([0, T])$.

Consider the following iterated Riemann--Stieltjes
integral

\vspace{-1mm}
\begin{equation}
\label{abcd1}
\int\limits_0^T
\int\limits_0^{s}
d{\bf f}_{\tau}^{(i_1)p}d{\bf f}_{s}^{(i_2)p},\ \ \ 
i_1,i_2=1,\ldots,m,
\end{equation}

\vspace{3mm}
\noindent
where $d{\bf f}_{\tau}^{(i)p}$ is defined by the
relation
(\ref{um1xxx}).

Let us substitute (\ref{um1xxx}) into (\ref{abcd1}) 

\vspace{-1mm}
\begin{equation}
\label{set18}
\int\limits_0^T
\int\limits_0^{s}
d{\bf f}_{\tau}^{(i_1)p}d{\bf f}_{s}^{(i_2)p}=
\sum\limits_{j_1,j_2=0}^p
C_{j_2 j_1} \zeta_{j_1}^{(i_1)}\zeta_{j_2}^{(i_2)},
\end{equation}

\vspace{3mm}
\noindent
where 
$$
C_{j_2 j_1}=
\int\limits_0^T \phi_{j_2}(s)\int\limits_0^s
\phi_{j_1}(\tau)d\tau ds
$$

\vspace{3mm}
\noindent
is the Fourier coefficient; another notations 
are the same as in (\ref{um1xxxx1}).

As we noted above, approximations of the Wiener process that are
similar to (\ref{um1xx})
were not considered in \cite{W-Z-1}, \cite{W-Z-2}
(also see Theorems 7.1, 7.2 in \cite{Watanabe}).
Furthermore, the extension of the results of Theorems 7.1 and 7.2
\cite{Watanabe} to the case under consideration is
not obvious.

On the other hand, we can apply Theorem 1 from this paper
and obtain from (\ref{set18}) the desired result

\vspace{-1mm}
\begin{equation}
\label{umen-bl}
\hbox{\vtop{\offinterlineskip\halign{
\hfil#\hfil\cr
{\rm l.i.m.}\cr
$\stackrel{}{{}_{p\to \infty}}$\cr
}} }
\int\limits_0^T
\int\limits_0^{s}
d{\bf f}_{\tau}^{(i_1)p}d{\bf f}_{s}^{(i_2)p}=
\hbox{\vtop{\offinterlineskip\halign{
\hfil#\hfil\cr
{\rm l.i.m.}\cr
$\stackrel{}{{}_{p\to \infty}}$\cr
}} }
\sum\limits_{j_1,j_2=0}^p
C_{j_2 j_1} \zeta_{j_1}^{(i_1)}\zeta_{j_2}^{(i_2)}=
{\int\limits_0^{*}}^T
{\int\limits_0^{*}}^s
d{\bf f}_{\tau}^{(i_1)}d{\bf f}_{s}^{(i_2)}.
\end{equation}

\end{document}